\documentclass[11pt]{amsart}
\usepackage{amssymb, amsmath, amsmath}
\usepackage{mathrsfs}
\usepackage{amscd}  
\usepackage{fullpage}
\usepackage[all]{xy} 
\usepackage{setspace}
\usepackage[dvipsnames]{xcolor}
\usepackage[margin=1.86  cm, headsep = 0.1 cm, footskip = 0.7 cm]{geometry}
\usepackage[english]{babel}
\usepackage{amsfonts}
\usepackage{xurl}
\usepackage[T1]{fontenc}
\usepackage{comment}
\usepackage[shortlabels]{enumitem}
\usepackage{tikz-cd}
\usepackage{quiver}
\usepackage{mathtools}
\usepackage{thmtools, thm-restate} 
\usepackage{hyperref}
\usepackage{csquotes}

\definecolor{darkred}{rgb}{0.5,0,0}
\definecolor{darkgreen}{rgb}{0,0.5,0}
\definecolor{darkblue}{rgb}{0,0,0.5}
\hypersetup{linkcolor=darkgreen,filecolor=darkgreen,urlcolor=darkgreen,citecolor=darkgreen,colorlinks}

\definecolor{c1}{RGB}{82,114,252}

\definecolor{c2}{RGB}{200,150,200}

\definecolor{c3}{RGB}{10,180,250}

\definecolor{c4}{RGB}{100,250,80}

\newcommand{\Z}{\mathbb{Z}}

\newcommand\C{\mathbb{C}}

\newcommand\Q{\mathbb{Q}}

\newcommand{\A}{\mathbb{A}} 
\newcommand\G{\mathbb{G}} 
\renewcommand{\H}{\mathrm{H}}

\newcommand{\rec}{\mathrm{rec}}

\newcommand{\et}{\text{\'et}} 
\newcommand{\lr}[1]{\langle #1 \rangle}

\DeclareMathOperator{\coker}{coker}

\DeclareMathOperator{\Spec}{Spec}

\DeclareMathOperator{\Br}{Br}

\DeclareMathOperator{\inv}{inv} 

\DeclareMathOperator{\Pic}{Pic}

\DeclareMathOperator{\Hom}{Hom}

\DeclareMathOperator{\im}{im} 

\DeclareMathOperator{\CH}{CH} 
\DeclareMathOperator{\Gal}{Gal}
 
\DeclareMathOperator{\cores}{cores} 
\DeclareMathOperator{\conn}{conn}

\DeclareMathOperator{\res}{res}

\newtheorem{theorem}{Theorem}[section]
\newtheorem{proposition}[theorem]{Proposition}
\newtheorem{lemma}[theorem]{Lemma}
\newtheorem{corollary}[theorem]{Corollary}
\theoremstyle{definition}
\newtheorem{definition}[theorem]{Definition}
\newtheorem{conjecture}[theorem]{Conjecture}
\theoremstyle{remark}
\newtheorem{remark}[theorem]{Remark}
\numberwithin{theorem}{section}
\numberwithin{equation}{section}

\usepackage{tcolorbox}
\usepackage{tikz-cd}
\usepackage{adjustbox}
\usepackage{marvosym}
\usepackage{ stmaryrd }
\usepackage{cancel}
\usepackage[normalem]{ulem}

\usepackage[OT2,T1]{fontenc}
\DeclareSymbolFont{cyrletters}{OT2}{wncyr}{m}{n}
\DeclareMathSymbol{\Sha}{\mathalpha}{cyrletters}{"58}
\DeclareMathSymbol{\locconstBr}{\mathalpha}{cyrletters}{"42}

\usepackage[
backend=biber,
style=alphabetic,
sorting=nyt, maxalphanames=99, minalphanames=99,
maxnames=50
]{biblatex}
\title{Connected descent and Brauer-Manin set for 0-cycles} 

\author{Francesca Balestrieri}
\address{Francesca Balestrieri\\ The American University of Paris\\
129 Rue de l'Universit\'{e}\\ 75007 Paris, France}
\email{fbalestrieri@aup.edu}


\let\oldtocsubsection\tocsubsection
\renewcommand{\tocsubsection}[3]{\hspace{1.5em}\oldtocsubsection{#1}{#2}{#3}}

\date{\today}

\begin{document}
\maketitle
\begin{abstract}
We answer a question by Hui Zhang and show that, for any smooth, quasi-projective, geometrically integral variety $X$ over a number field $k$, the connected descent set for 0-cycles is equal to the Brauer-Manin set for 0-cycles. In particular, the connected descent and the abelianised connected descent sets for 0-cycles coincide. In the same spirit, we also show that the descent set under all tori for 0-cycles is equal to the algebraic Brauer-Manin set for 0-cycles. Finally, we give a new formulation of Conjecture (E) using the connected descent sets. 
\end{abstract}

\section{Introduction}
The motivation for this paper comes from trying to understand whether some facts that hold true for rational points also hold for 0-cycles. Specifically, for rational points there are many equivalence results between "descent-type" obstruction sets (i.e. obstruction sets using descent under certain classes of linear algebraic groups) and "Brauer-type" obstruction set (i.e. obstruction sets using subgroups of the Brauer group). For example, if $X$ is a smooth, quasi-projective, geometrically integral variety over a number field $k$, we have that the connected descent set 
\[X(\A_k)^{\textrm{conn}} := \bigcap_{G \textrm{ connected lin alg $k$-group}} \bigcap_{[f: Y \to X] \in H^1(X, G)} \bigcup_{\sigma \in H^1(k, G)} f^\sigma(Y^\sigma(\A_k))\]
is equal to the Brauer-Manin set $X(\A_k)^{\Br}$ by a result by Harari \cite[Thm 2(ii)]{Harari02}; similarly, we have that the descent set under all linear algebraic groups of multiplicative type over $k$ is the same as the algebraic Brauer-Manin set $X(\A_k)^{\Br_1}$, where $\Br_1(X):= \ker(\res_{\overline{k}/k}: \Br X \to \Br X_{\overline{k}})$, see \cite{Sko-bielliptic-99} \cite{Harari_Skorobogatov_2013}; for other examples, see e.g. \cite{Dem09}, \cite{Sko-Descent-EtBrauer}, \cite{CDX19}, \cite{Balestrieri16}. For 0-cycles, the formalism for a theory of descent obstructions has been recently developed in \cite{balestrieri-berg}, \cite{Zhang-0cycles-25}, \cite{BGNSST26}. In particular, in \cite{Zhang-0cycles-25}, using his abelianised descent set, Zhang proved the following analogue for 0-cycles of the result by Harari for rational points.
\begin{theorem}[{\cite[Theorem 4.7]{Zhang-0cycles-25}}]\label{thmintro1}
Let $X$ be a smooth and geometrically integral variety over a number field $k$, then
\(Z_{0, \A} (X)^{\Br} = Z_{0, \A} (X)_{ab}^{\conn}, \)
where $ Z_{0, \A} (X)_{ab}^{\conn}$ is the abelianised connected descent set (see \cite{Zhang-0cycles-25} for more details).     
\end{theorem}
For $X$ projective, in \cite[Remark 5.8]{Zhang-0cycles-25} Zhang also asked under which reasonable assumptions we have $\overline{Z_{0, \Omega}(X)^{\conn}} = Z_{0, \Omega}(X)^{\Br}$, where $\overline{Z_{0, \Omega}(X)^{\conn}}$ denotes (a variant of ) the Balestrieri-Berg connected descent set (see Definition \ref{mydef2}) and where the closure is with respect to the weak approximation topology (see e.g. \cite{Zhang-0cycles-25} for a definition). In \cite{BGNSST26}, using their refined descent formalism, the authors answer Zhang's question under some natural conditions.

\begin{theorem}[{\cite[Theorem 1.8]{BGNSST26}}]
Let $X$ be a smooth, projective, geometrically integral variety over a number field $k$. Suppose that
\begin{enumerate}[(i)]
\item \label{part:conndense-cond1} $Z_{0, \Omega}(X)^{\Br}$ is open in $Z_{0, \Omega}(X)$ with respect to the weak approximation topology, and 
\item \label{part:conndense-cond2} for every positive integer $d$, there exists a finite extension $k'/k$ such that for any degree $d$ extension $K/k$ which is linearly disjoint from $k'/k$, the restriction map
\[
\Br(X)/\Br_0(X)\ \to\ \Br(X_K)/\Br_0(X_K)
\]
is surjective.
\end{enumerate}
Then $\overline{Z_{0, \Omega}(X)^{\conn}} = Z_{0, \Omega}(X)^{\Br}$ and $\overline{Z_{0, \Omega}^d (X)^{\conn}} = Z_{0, \Omega}^d (X)^{\Br}$ for all $d\in\Z$.
\end{theorem}

The main aim of this paper is to answer Zhang's question in a much stronger and more general way, by removing all conditions and by showing that we don't need to take the closure. Our main result is the following.

\begin{theorem}[Corollary \ref{maincor}]\label{thmintro2}
     Let $X$ be a smooth, quasi-projective, geometrically integral variety over a number field $k$. Let $d \in \Z$. Then  $Z_{0, \A}^d(X)^{\conn} = Z_{0, \A}^d(X)^{\Br}$ and $Z_{0, \A}^d(X)^{\text{\normalfont tori}} = Z_{0, \A}^d(X)^{\Br_1}$.
\end{theorem}

In particular, combining Theorems \ref{thmintro1} and \ref{thmintro2}, we have $Z_{0, \A}^d(X)^{\conn} =Z_{0, \A}^d(X)_{ab}^{\conn}$ for any $d \in \Z$, i.e. the (variant of the) Balestrieri-Berg connected descent set defined in Definition \ref{mydef2} is the same as Zhang's abelianised connected descent set for any smooth, quasi-projective, geometrically integral variety $X$ over $k$. In fact, a stronger result is true.
\begin{theorem}[Corollary \ref{cor:zhang}]
  Let $X$ be a smooth, quasi-projective, geometrically integral variety over a number field $k$. Let $G$ be a connected linear algebraic group over $k$ and let $f: Y \to X$ be a torsor under $G$. Then, for any $d \in \Z$, we have
  \[ Z_{0, \A}^d(X)^f =  Z_{0, \A}^d(X)^f_{ab},\]
  where $ Z_{0, \A}^d(X)^f_{ab}$ is the abelianised descent set as in \cite{Zhang-0cycles-25}.
\end{theorem}

The second part of the paper is about reformulating one of the main conjectures on the qualitative arithmetic of 0-cycles using our main result above.
Recall that the \emph{Chow group} $\CH_0(X)$ is the quotient of $Z_0(X)$ by the subgroup of $0$-cycles rationally equivalent to zero. For 0-cycles on smooth, projective, geometrically integral varieties over number fields, there is a main conjecture on their qualitative arithmetic behaviour (originally due to Colliot-Thélène, Sansuc, Saito and Kato, and reformulated in its current form by van Hamel and Wittenberg, see \cite{CTS81}, \cite{KatoSaito86}, \cite{CT-Chow-93}, \cite{vanHamel03}, \cite{Wit12}) stating essentially that the Brauer-Manin obstruction completely captures their behaviour, in stark contrast to what happens for rational points, where there are numerous counterexamples - see e.g. \cite{Sko-bielliptic-99}, \cite{Poonen10} - where the Brauer-Manin obstruction is not sufficient. 

\begin{conjecture}[Conjecture (E)] \label{conjE}   For any smooth, projective, geometrically integral variety $X$ over a number field $k$, the following complex is exact,
\[ \widehat{\CH_0(X)} \longrightarrow \prod_{v \in \Omega_k} \widehat{\CH_0'(X_{k_v})} \longrightarrow \Hom(\Br X, \Q/\Z),\]
where $\widehat{\ M\ }:= \varprojlim_{n} M/n M$, and where $\CH_0'(X_{k_v})$ coincides with $\CH_0(X_{k_v})$ at finite places and at infinite places is the modification $\CH_0'(X_{k_v}) := \coker(N_{\overline{k}_v/k_v})$, where $N_{\overline{k}_v/k_v}: \CH_0(X_{\overline{k}_v}) \to \CH_0(X_{k_v})$ is the norm map (see \cite{Wit12} or \cite{CT-Chow-93} for further details).
\end{conjecture}
  In the light of our main theorem, it makes sense to ask if we can recast Conjecture (E) in terms of the connected descent. We prove the following result. 

  \begin{theorem}[Theorem \ref{conjE}] Let $X$ be a smooth, projective, geometrically integral variety over a number field $k$. 
Let $\phi: Z_{0, \Omega}(X) \to \widehat{\prod_v \CH'_0(X_v)}$ be given by $(z_v)_v \mapsto (([z_{v}]_n)_{n \geq 1})_v$, where $[- ]_n$ denotes the class in $\CH'_0(X_{k_v})/n$.
    Then Conjecture (E) is equivalent to 
    \[
    \im\left(\widehat{\CH_0(X)} \to  \widehat{\prod_v \CH'_0(X_v)} \right) = \bigcap_{G \textrm{ connected lin alg $k$-groups}} \bigcap_{[f: Y \to X] \in H^1(X, G)} \overline{\phi(Z_{0, \Omega}(X)^f)}.
    \]
    In particular, it suffices to consider the split groups $G = PGL_{n, k}$ for $n \geq 1$. 
\end{theorem}

\section{Preliminaries}

Throughout this section, we let $k$ be a number field and let $X$ be a smooth, quasi-projective, geometrically integral variety defined over $k$. 
We write $X_0$ for the set of closed points on $X$. 
We write $X(\A_k)$ for the set of adelic points on $X$, and we let $X(k_\Omega) \coloneq\Pi_{v \in \Omega} X(k_v)$.
Then $X(\mathbb{A}_k) \subset X(k_\Omega)$, and equality holds when $X$ is proper.

\begin{definition}
A \emph{0-cycle $z$ on $X$} is a formal sum
$z \coloneq\sum_{\substack{x\in X_0}} n_x x,$
with $n_x\in\mathbb{Z}$, such that $n_x = 0$ for all but finitely many $x \in X_0$. The degree of $z$ is given by
$\deg(z) \coloneq\sum_{\substack{x\in X_0
}} n_x [k(x):k].$

\end{definition}
We write $Z_0(X)$ for the \emph{set of 0-cycles on $X$} and $Z_0^d(X)$ for the \emph{subset of 0-cycles of degree $d$}. 

\begin{remark} 
The set $Z_0(X)$ is an abelian group under addition. For 0-cycles $z_1,z_2$ on $X$, we have the equality $\deg(z_1+z_2) = \deg(z_1) + \deg(z_2)$. 
\end{remark}

\begin{definition} 
A 0-cycle $z = \sum_{x\in X_0} n_x x$ on $X$ is called \emph{effective} if $n_x\geq 0$  for all $x\in X_0$. We say that $x \in X_0$ is in the \emph{support of $z$} if $n_x \neq 0$.
\end{definition}

\begin{definition}
We write $Z^{d}_{0,\Omega}(X) \coloneq\prod_{v \in \Omega} Z^{d}_0(X_{k_v})$ and $Z_{0,\Omega}(X) \coloneq\bigsqcup_{d \in \Z} Z^{d}_{0,\Omega}(X)$. 
We call elements of $Z_{0,\Omega}(X)$ \emph{$\Omega$-local 0-cycles} on $X$.
\end{definition}

\begin{definition}
    An \emph{adelic 0-cycle of degree $d$} is a local 0-cycle $(z_v)_v \in Z_{0, \Omega}^d(X)$ satisfying the following conditions: there is a finite set of places $S \subset \Omega_k$ containing the infinite places and a separated model $\mathcal X$ of finite type over $\mathcal O_{k, S}$ such that, for any $v \notin S$ and every $x_v \in supp(z_v)$, we have that the morphism $\Spec k_v(x_v) \to X$ extends to a morphism $\Spec O_{k_v(x_v)} \to \mathcal X$. We write $Z_{0, \A}^d(X)$ for the set of adelic 0-cycles of degree $d$ on $X$, and write $Z_{0, \A}(X) := \bigsqcup_{d \in \Z} Z_{0, \A}^d(X)$. 
\end{definition}

\begin{remark}
    If $X$ is proper, then $Z^d_{0, \A}(X) = Z^d_{0, \Omega}(X)$ for all $d \in \Z$.
\end{remark}

\begin{remark}
For any degree $d$, we have a natural inclusion \(
Z_0^d(X) \subset Z_{0,\A}^d(X) \subset Z_{0,\Omega}^d(X).\)
\end{remark}

\subsection{Brauer--Manin sets for rational points and 0-cycles} 
\label{subsection:obstruction-sets}
Let \(X\) be a smooth, quasi-projective, geometrically integral variety over $k$.
The \emph{Brauer group of \(X\)} is defined as \(\Br (X) \coloneq\H^2_{\et}(X, \G_m)\). The \emph{algebraic Brauer group of $X$} is $\Br_1 X := \ker(\res_{\overline{k}/k}: \Br X \to \Br (X_{\overline{k}}))$. 
Let $T$ be a $k$-scheme with $X(T)\neq\emptyset$. Viewing $x\in X(T)$ as a morphism $x:T\to X$ induces a pullback map $x^*:\Br(X)\to\Br(T)$. An element $\alpha\in \Br(X)$ induces an evaluation map $\alpha:X(T)\to\Br(T)$ given by $x\mapsto \alpha(x):=x^*(\alpha)$.
We have the commutative diagram
\begin{center}
\begin{tikzcd}
\label{fig:Brauer-rational-points}
	& {X(k)} & X(\A_k) \\
	0 & {\Br(k)} & {\bigoplus_{v\in\Omega} \Br(k_v)} && {\mathbb{Q}/\mathbb{Z}} & 0, 
	\arrow[hook, from=1-2, to=1-3]
	\arrow["\alpha", from=1-2, to=2-2]
	\arrow["{\prod_{v\in\Omega} \alpha}", from=1-3, to=2-3]
	\arrow[from=2-1, to=2-2]
	\arrow[from=2-2, to=2-3]
	\arrow["{ \sum_{v\in\Omega} \text{inv}_v}", from=2-3, to=2-5]
	\arrow[from=2-5, to=2-6]
\end{tikzcd}
\end{center}
where the sequence at the bottom is the fundamental short exact sequence coming from class field theory (see e.g. \cite[Thm 12.1.8(iii)]{CTS-BrauerBook}). In \cite{Manin71}, Manin defined the now called \emph{Brauer--Manin pairing} $
\lr{ \cdot, \cdot}_{BM}: X(\A_k) \times \Br(X) \to \Q/\Z$ given by $((x_v)_v, \alpha ) \mapsto \sum_{v \in \Omega} \mathrm{inv}_v(\alpha(x_v))$. 
For $\alpha \in \Br(X)$, we define the \emph{set of adelic points orthogonal to $\alpha$} to be
\[X(\A_k)^{\alpha}\coloneq\{x\in X(\A_k)\;:\; \lr{x,\alpha}_{BM}=0\}.\]
For a subset $B \subset \Br(X)$, we define the \emph{set of adelic points orthogonal to $B$} to be
\[X(\A_k)^{B }=\bigcap_{\alpha\in B} X(\A_k)^{\alpha}.\]

Now let \(x_v \in (X_{k_v})_0\) be a closed point with residue field $k_v(x_v)$.  
For \(\alpha \in \Br(X)\), evaluation at $x_v$ yields
$\alpha(x_v) \in \Br(k_v(x_v))$, and then corestriction $\mathrm{cores}_{k_v(x_v)/k_v} \colon \Br(k_v(x_v)) \to \Br(k_v)$ yields an element in $\Br(k_v)$.
If \(z_v = \sum n_{x_v} x_v\) is a 0-cycle, then we can apply this process to each \(x_v\), 
allowing us to extend the Brauer--Manin pairing to a pairing $\lr{ \cdot, \cdot}_{BM}: Z_{0, \A}(X) \times \Br (X) \to \Q/\Z$ given by 
\[\Biggl(\Bigl(\sum_{\substack{x_v \in (X_{k_v})_0}}  n_{x_v} x_v\Bigr)_v\ ,\ \alpha \Biggr)  \mapsto \sum_{v\in \Omega} \sum_{\substack{x_v \in (X_{k_v})_0
}}  n_{x_v} \text{inv}_v(\mathrm{cores}_{k_v(x_v)/k_v}(\alpha(x_v))).\]

\begin{definition}[{\cite{CT-Chow-93}}]
\label{defn:unrefined-Brauer}  
For $\alpha\in \Br(X)$, define the \emph{group of adelic 0-cycles on $X$ orthogonal to $\alpha$} to be
$Z_{0, \A}(X)^\alpha \coloneq\left\{ z \in Z_{0, \A}(X)\;:\; \lr{z,\alpha}_{BM} = 0 \right\}.$
For a subset $B \subset \Br(X)$, we define the \emph{group of adelic 0-cycles on $X$ orthogonal to $B$} to be $ Z_{0, \A}(X)^{B} \coloneq\bigcap_{\alpha \in B} Z_{0, \A}(X)^\alpha.$
For each $d\in\Z$, we define $Z_{0, \A}^d(X)^{B}\coloneq Z_{0, \A}(X)^{B} \cap Z_{0, \A}^d(X)$.
\end{definition}

\subsection{Descent sets and recombining maps}\label{sec:descent}
Let $Y$ be a quasi-projective $k$-variety endowed with an action of an algebraic group $G$ over $k$. Then we denote by $Y^\sigma$ the twist of $Y$ by $\sigma \in \H^1(k,G)$, see \cite[p.12]{Sko-Torsors}.

\begin{definition}\label{defn:f-descent_pt}
Let $G$ be a linear algebraic group over $k$ and let $f: Y\rightarrow X$ be an $X$-torsor under $G$ defined over $k$~\cite[Definition {2.2.1}]{Sko-Torsors}. 
For any $\sigma\in \H^1(k,G)$, we have an $X$-torsor under $G^\sigma$ of the form $f^\sigma : Y^\sigma \rightarrow X$.
We can partition the set of $k$-rational points on $X$ as 
\begin{equation}\label{eq:partition}
X(k) = \bigsqcup_{\sigma\in \H^1(k,G)} f^\sigma(Y^\sigma(k)),
\end{equation}
see e.g.\ \cite[\S 2]{Sko-Torsors}.
The \emph{$f$-descent set for rational points} is defined as
\begin{equation*}
 X(\A_k)^f \coloneq \bigcup_{\sigma\in \H^1(k,G)} f^\sigma(Y^\sigma(\A_k)) = \{ (x_v)_v: \exists \sigma \in H^1(k, G) \text{ such that } [Y_{x_v}] = \sigma \ \forall v \in \Omega_k\},
\end{equation*}
where $Y_{x_v} = f^{-1}(\{ x_v\})$ denotes the fibre above $x_v$.
\end{definition}

\begin{remark}
The inclusions $Y^\sigma(k) \subset Y^\sigma (\A_k)$, together with~\eqref{eq:partition}, yield an inclusion $X(k) \subset X(\A_k)^f$.
\end{remark}

To define an analogue of partition~\eqref{eq:partition} for 0-cycles, we first introduce the global and local recombining maps defined in \cite{balestrieri-berg}.

\begin{definition}
Let $G$ be a linear algebraic group over $k$, 
let $f:Y\rightarrow X$ be an $X$-torsor under $G$ defined over $k$ and let $L$ be a finite extension of $k$. 
For any $\sigma\in \H^1(L,G)$, we have an $X_L$-torsor under $G_L^\sigma$ defined over $L$ given by $f_L^\sigma:Y_L^\sigma\rightarrow X_L$. This twisted $X_L$-torsor induces pushforward maps 
$f_{L,*}^\sigma : Z_0(Y_L^\sigma) \rightarrow Z_0(X_L)$ and $f_{L,*}^\sigma : Z_{0, \A}(Y_L^\sigma) \rightarrow Z_{0, \A}(X_L)$ 
which preserve degrees of 0-cycles. 
\end{definition}

\begin{definition}
Let $L$ be a finite extension of $k$. There exists a natural \emph{corestriction map}
$$N_{L/k}: Z_0(X_L) \rightarrow Z_0(X),$$
which sends a 0-cycle of degree $d$ on $X_L$ to a 0-cycle of degree $[L:k]\cdot d$ on $X$.
The local corestriction maps $N_{L_w/k_v}: Z_0(X_{L_w}) \rightarrow Z_0(X_{k_v})$ for $w\in\Omega_L$ lying above $v\in\Omega_k$ yield a map 
\[
N_{L/k}^\A : Z_{0, \A}(X_L) \rightarrow Z_{0, \A}(X)
\]
given by $(z_w)_{w\in\Omega_L} \mapsto \bigl(\sum_{w \mid v} N_{L_w/k_v}(z_w)\bigr)_{v\in\Omega_k}$ 
that restricts to $N_{L/k}$ when one views $Z_0(X_L)$ and $Z_0(X)$ as subgroups of $Z_{0, \A}(X_L)$ and $Z_{0, \A}(X)$, respectively. Explicitly, for $z\in Z_0(X_L)$, we have that 
$
N_{L/k}^\A((z)_{w\in\Omega_L})=(N_{L/k}(z))_{v\in\Omega_k}.$
\end{definition}

\begin{definition}
Let $G$ be a linear algebraic group over $k$ and let $f: Y\rightarrow X$ be an $X$-torsor under $G$ defined over $k$. Define the \emph{global $f$-recombining map} $$ \rec_f : \bigoplus_{L/k \text{ finite}} \bigoplus_{\sigma\in\H^1(L,G)} Z_0(Y_L^\sigma) \rightarrow Z_0(X)$$
in the following way. Let $L_1,\dots, L_n$ be a finite collection of finite extensions of $k$. For each $L_i$ let $\sigma_{i,1},\dots ,\sigma_{i, m_i}$ be a finite collection of elements of $\H^1(L_i,G)$. For a 0-cycle $z\in Z_0(Y_{L_i}^{\sigma_{i,j}})$, the pushforward $f_{L_i,*}^{\sigma_{i,j}}(z)$ lies in $Z_0(X_{L_i})$. We
define $\rec_f(z)\coloneq N_{L_i/k}(f_{L_i ,*}^{\sigma_{i,j}}(z))\in Z_0(X)$ and extend by linearity.
\end{definition}

\begin{remark}
Up to a change of notation, this is the map $g_*$ described in Definition 2.2 of \cite{balestrieri-berg}.
\end{remark}

Using this recombining map, we define the analogue of partition~\eqref{eq:partition} for 0-cycles as follows.

\begin{lemma}\label{lem:partition0cycles}
Let $G$ be a linear algebraic group over $k$ and let $f: Y\rightarrow X$ be an $X$-torsor under $G$ defined over $k$. Then 
\begin{equation*}
Z_0(X) = \rec_f \left( \bigoplus_{\substack{L/k\ 
\mathrm{finite}}}  \bigoplus_{\sigma \in \H^1(L,G)} Z_0(Y_L^\sigma)\right).
\end{equation*}
\end{lemma}
\begin{proof}
See~\cite[Proposition~2.5]{balestrieri-berg}.
\end{proof}

To define descent sets for 0-cycles, we also require the notion of a local recombining map.

\begin{definition}[{\cite[Def 3.1]{balestrieri-berg}}]
Let $G$ be a linear algebraic group over $k$ and let $f: Y\rightarrow X$ be an $X$-torsor under $G$ defined over $k$. 
Define the \emph{local $f$-recombining map}
$$ \rec_f^\A : \bigoplus_{L/k \text{ finite}} \bigoplus_{\sigma\in\H^1(L,G)} Z_{0,\A}(Y_L^\sigma) \rightarrow Z_{0,\A}(X)$$
in the following way. Let $L_1,\dots, L_n$ be a finite collection of finite extensions of $k$. For each $L_i$ let $\sigma_{i,1},\dots ,\sigma_{i, m_i}$ be a finite collection of elements of $\H^1(L_i,G)$. For an adelic 0-cycle $(z_w)_{w\in\Omega_{L_i}}\in Z_{0,\A}(Y_{L_i}^{\sigma_{i,j}})$, define $\rec_f^\A((z_w)_{w\in\Omega_{L_i}})\coloneq N_{L_i/k}^\Omega(f_{L_i ,*}^{\sigma_{i,j}}((z_w)_w))\in Z_{0,\A}(X)$ and extend by linearity. 
\end{definition}

\begin{definition}[{\cite[Definition 3.3]{balestrieri-berg}}]
\label{defn:non-refined-f-descent}
Let $G$ be a linear algebraic group over $k$ and let $f:Y\rightarrow X$ be an $X$-torsor under $G$ defined over $k$. The \emph{$f$-descent set for 0-cycles} is defined as 
\begin{equation*}
Z_{0,\A}(X)^f \coloneq\rec^\A_f \left( \bigoplus_{L/k \textrm{ finite}}  \bigoplus_{\substack{\sigma\in \H^1(L,G)}} Z_{0,\A}(Y_L^\sigma)\right).
\end{equation*}
\end{definition}
\begin{remark}
    Lemma~\ref{lem:partition0cycles} shows that \(Z_0(X)\subset Z_{0,\A}(X)^f.\)
\end{remark}

\begin{definition} \label{mydef2} If $\mathcal{G}$ is a class of linear algebraic groups over $k$, we let
\[
Z_{0,\A}(X)^{\mathcal{G}} := \bigcap_{G \in \mathcal{G}} \bigcap_{[f: Y \to X] \in H^1(X, G)} Z_{0,\A}(X)^f.
\]
\end{definition}

\begin{remark}
We note that Definition \ref{mydef2} is slightly different from the one in \cite{balestrieri-berg}, because of a different order of the quantifiers involved: in \cite{balestrieri-berg}, they first consider a finite set of finite extensions of $k$, and then that finite set of extensions needs to work for all linear algebraic groups and all torsors. So the descent set in \cite{balestrieri-berg} is potentially smaller than the set defined in Definition \ref{mydef2}.
\end{remark}


\section{Torsors under connected linear algebraic groups
}

Let $X$ be a smooth, quasi-projective, and geometrically integral variety over a number field $k$. Let $G$ be a connected linear algebraic group over $k$ and let $f: Y \to X$ be a torsor under $G$. Central extensions of $G$ by $\G_m$ given by
\[
1 \to \G_{m, k} \to E \to G \to 1,
\]
for some linear algebraic group $E$ over $k$ are classified by the abelian group $Ext_{k-groups}^c(G, \G_{m, k})$.
Since $G$ is connected, we have an isomorphism from $Ext_{k-groups}^c(G, \G_{m, k}) \xrightarrow{\sim} \Pic(G)$ (see e.g. \cite[Prop 2.9]{CTXu09}). Let $u \in \Pic(G)$ and let $\partial_u: H^1(X, G) \to H^2(X, \G_m) = \Br X$ be the coboundary map associated with the central extension corresponding to $u$ via the isomorphism $Ext_{k-groups}^c(G, \G_{m, k}) \xrightarrow{\sim} \Pic(G)$.
We define the map 
\[ \Delta_{Y/X} : \Pic(G) \cong Ext_{k-groups}^c(G, \G_{m, k})  \to \Br X\]
given by $u \mapsto \partial_u([f: Y \to X])$. Since $\Delta_{Y/X}$ is compatible with base change, we have in particular that for any $x \in X$, the evaluation $ev_x(\partial_u([f: Y \to X]))$ is the same as $\partial_u([f_x : Y_x \to \{x\}]) \in \Br k(x)$.

\begin{lemma} \label{lem:lift} Let $X$ be a smooth, quasi-projective, and geometrically integral variety over a number field $k$. Let $G$ be a connected linear algebraic group over $k$ and let $f: Y \to X$ be a torsor under $G$.  Let $L/k$ be a finite extension. Then $(x_w) \in X_L(\A_L)$ lifts to an adelic point on $Y^\sigma_L(\A_L)$ for some $\sigma \in H^1(L, G)$ if and only if $(x_w)_w$ is orthogonal to $\Delta_{Y_L/X_L}(\Pic(G_L)) \subset \Br(X_L)$.
\end{lemma}
\begin{proof}
   By the exactness of the sequence in equation (3.8) and by the diagram in equation (3.9) in the proof of \cite[Thm 3.5]{CDX19}, we have that $\ker(f_L^\ast: \Br(X_L) \to \Br(Y_L)) = \Delta_{Y_L/X_L}(\Pic(G_L))$ and that
   \[ X_L(\A_L)^{\ker(f_L^\ast)} = \bigcup_{\sigma \in H^1(L,G)} f_L^\sigma(Y^\sigma_L(\A_L)),\]
   as required.
\end{proof}

\begin{lemma} \label{lem0} Let $X$ be a smooth, quasi-projective, and geometrically integral variety over a number field $k$. Let $G$ be a connected linear algebraic group over $k$ and let $f: Y \to X$ be a torsor under $G$. Then $\Pic(G)$ is finite and there is a finite Galois extension $F/k$ such that, for any finite extension $L/k$ linearly disjoint from $F$ over $k$, we have an isomorphism $\res_{L/k}: \Pic(G) \xrightarrow{\sim} \Pic(G_L)$ and $\Delta_{Y_L/X_L}(\Pic(G_L)) = \res_{L/k}(\Delta_{Y/X}(\Pic(G)))$.    
\end{lemma}
\begin{proof}
    Since $k$ is perfect and $G$ is connected and affine, by \cite[Lemma 6.9(i)]{Sansuc81} $\Pic(G)$ is finite. By \cite[Prop-déf 3.1]{CTFlasque08}  there is a flasque resolution of $G$ given by
    \[
    1 \to T \to H \to G \to 1, 
    \]
    where $T$ is a flasque torus and $H$ is quasi-trivial. Let $P:=H^{tor}$ denote the maximal torus quotient of $H$, and let $\alpha : T \to H \to P$. Then by \cite[Prop 3.3]{CTFlasque08}, for any $L/k$ we get the exact sequence compatible with base change to any extension of $k$ given by
    \begin{equation}\label{flasque}
    \widehat{P}^{\Gamma_L} \xrightarrow{\alpha^\ast} \widehat{T}^{\Gamma_L} \xrightarrow{d_L} \Pic(G_L) \to 0,
    \end{equation}
    with the property that $d_L(\chi_L) = \res_{L/k}(d_k(\chi))$ for any $\chi \in \widehat{T}^{\Gamma_k}$.
    Now let $F/k$ be a Galois finite extension splitting both $P$ and $T$. For any $L/k$ linearly disjoint from $F$ over $k$, we have that the composite $\Gamma_L \to \Gal(FL/L) \xrightarrow{\sim} \Gal(F/k)$ is surjective. Since the Galois actions on both $\widehat{P}$ and $\widehat{T}$ factor by construction through $\Gal(F/k)$, we have that $\widehat{P}^{\Gamma_L} = \widehat{P}^{\Gamma_k}$ and $\widehat{T}^{\Gamma_L} = \widehat{T}^{\Gamma_k}$. Hence, by \eqref{flasque} and by the compatibility of the map $d_k$ with base change $\res_{L/k}$, we have that the natural restriction map
    \[
    \res_{L/k}: \frac{\widehat{T}^{\Gamma_k}}{\alpha^\ast(\widehat{P}^{\Gamma_k})} = \Pic(G_k) \to \frac{\widehat{T}^{\Gamma_L}}{\alpha^\ast(\widehat{P}^{\Gamma_L})} = \Pic(G_L)
    \]
    is an isomorphism. 
    For $u \in \Pic(G)$, let \[
    1 \to \G_m \to E \to G \to 1
    \]
be the corresponding central extension. Let $s_L: X_L \to X$. Then $\Delta_{Y_L/X_L}(u_L) = \partial_{u_L}([f_L: Y_L \to X_L]) = s_L^\ast\partial_u([f: Y \to X]) = s_L^\ast(\Delta_{Y/X}(u))$. As we vary $u \in \Pic(G)$, by the isomorphism $\res_{L/k}: \Pic(G) \xrightarrow{\sim} \Pic(G_L)$ we also have that $u_L$ varies in $\Pic(G_L)$. Taking the images under $\Delta_{Y_L/X_L}$ yields $\Delta_{Y_L/X_L}(\Pic(G_L)) = s_L^\ast (\Delta_{Y/X}(\Pic(G))) = \res_{L/k}(\Delta_{Y/X}(\Pic(G)))$, as required.  
\end{proof}

\begin{lemma} \label{lem1}
    Let $X$ be a quasi-projective variety over $k$. Let $F/k$ be a finite Galois extension, and let $S \subset \Omega_k$ be a finite set of places of $k$. Let $(z_v)_{v \in S}$ be effective 0-cycles of a same degree $d > 0$. Then there is an extension $L/k$ linearly disjoint from $F$ over $k$ with $[L:k] = d$ and local points $x_w \in X(L_w)$ for each $w|v$ with $v \in S$ which recombine to $z_v$. If moreover we have local points $u_v \in X(k_v) \neq \emptyset$ for $v \notin S$ which are integral almost everywhere on a fixed model, then setting $x_w:=(u_v)_{L_w}$ for $w|v$ and $v \notin S$ yields an adelic point in $X(\A_L)$.
\end{lemma}
\begin{proof}
    Write $z_v = \sum_{u_v} n_{u_v} u_v$, and consider the finite étale algebra over $k_v$ of dimension $d$ given by $K_v := \prod_{u_v \in supp(z_v)} \prod_{i=1}^{n_{u_v}} k_v(u_v)$. Since $k_v$ is infinite, for each $(u_v, i)$ we can choose a primitive element of $k_v(u_v)/k_v$ so that their monic minimal polynomials $m_{u_v, i}$ are all distinct. Let $m_v := \prod_{u_v, i} m_{u_v, i}$, so that $m_v$ is separable and $\deg(m_v) = d$. Then we have that $K_v = k_v[t]/(m_v)$ and the divisor $D_v:= V(m_v)$ is reduced of degree $d$.
    By \cite[Prop 3.3.3]{Liang-Arithmetic-13} applied to $D_v$ and to the Hilbert subset defined by $\mathbb{A}^1_F \to \mathbb{A}^1_k$, we get a closed point $\theta \in \mathbb{A}^1_k$ such that, if $L:= k(\theta)$, we have that $[L:k] = d$ and $L$ is linearly disjoint from $F$ over $k$; moreover, for each $v \in S$ we have by Krasner's lemma that $\prod_{w|v} L_w = L \otimes_k k_v = K_v = \prod_{u_v} \prod_{i=1}^{n_{u_v}} k_v(u_v)$. Indexing the places $w$ above $v$ using the pairs $(u_v, i)$, say $w := w_{u_v, i}$, we fix the $k_v$-isomorphisms $\phi_{u_v, i}: k_v(u_v) \to L_{w_{u_v, i}}$ and we define local points $x_{w_{u_v, i}} \in X(L_{w_{u_v, i}})$ by 
    \[ \Spec L_{w_{u_v, i}} \xrightarrow{\phi_{u_v, i}^\ast} \Spec k_v(u_v) \xrightarrow{u_v} X_{k_v}. \]
    Then the pushforward under the base change $s_{L_w/k_v} : X_{L_w} \to X_{k_v}$ is 
    \[
    \sum_{w|v} s_{L_w/k_v, \ast}(x_w) = \sum_{u_v \in supp(z_v)}\sum_{i=1}^{n_{u_v}} [L_{w_{u_v, i}}: \phi_{u_v, i}(k_v(u_v))] u_v = \sum_{u_v} n_{u_v} u_v = z_v,  
    \]
    so the $(x_w)_{w|v}$ recombine to $z_v$.
    Finally, the last statement of the lemma is clear.  
\end{proof}

\begin{theorem} \label{thm:main} Let $X$ be a smooth, quasi-projective, geometrically integral variety over a number field $k$. Let $G$ be a connected linear algebraic group over $k$ and let $f: Y \to X$ be a torsor under $G$. Let $d \in \Z$. Then
\[
Z_{0,\A}^d(X)^f = Z_{0,\A}^d(X)^{\Delta_{Y/X}(\Pic(G))}.
\]
In particular, $Z_{0,\A}^d(X)^{H^1(X, G)} = Z_{0,\A}^d(X)^{\langle \Delta_{Y/X}(\Pic(G)) : [f: Y \to X] \in H^1(X, G) \rangle}$ and, if $\mathcal{G}$ is a class of connected linear algebraic groups over $k$, then 
\[
Z_{0,\A}^d(X)^{\mathcal{G}} = Z_{0,\A}^d(X)^{\langle \Delta_{Y/X}(\Pic(G)) : [f: Y \to X] \in H^1(X, G), G \in \mathcal{G} \rangle}.
\]
    
\end{theorem}
\begin{proof} We first prove the inclusion $\subseteq$. By additivity, it suffices to prove the inclusion for $N^\A_{L/k}((f^\sigma_L)_\ast z)$ for 
any finite extension $L/k$, $\sigma \in H^1(L, G)$, and $z = (z_w)_w \in Z_{0,\A}^e(Y_L^\sigma)$, for any $e \in \Z$. We fix such $L/k$, $\sigma \in H^1(L, G)$, $e \in \Z$, and  $z \in Z_{0,\A}^e(Y_L^\sigma)$.
Let $u \in \Pic(G)$ and let $\alpha := \Delta_{Y/X}(u) = \partial_u([f: Y \to X]) \in \Br X$. Say $\sigma = [g_\sigma: Z_\sigma \to \Spec L]$. Then we have that $(f_L^\sigma)^\ast Y_L = s^\ast(Z_\sigma)$, where $s: Y^\sigma_L \to \Spec L$ is the structure morphism, and that $(f^\sigma_L)^\ast \alpha_L = s^\ast (\partial_{u_L}(\sigma)) \in \Br(Y^\sigma_L)$, where $\beta:= \partial_{u_L}(\sigma) \in \Br(L)$ (see e.g. \cite[(3.7)-(3.9)]{CDX19} and \cite[p.314]{CTXu09}). Write $z_w = \sum_{y_w} n_{y_w} y_w$. Then for any $v \in \Omega_k$ and any $w | v$, we have
\[
\inv_v \alpha_{k_v}\left(N_{L_w/k_v}((f^\sigma_{L_w})_\ast z_w) \right) = \inv_w (((f_{L_w}^\sigma)^\ast \alpha_{L_w})(z_w)) = \sum_{y_w} n_{y_w} [L_w(y_w): L_w] \inv_w(\beta_{L_w}) = e \inv_w(\beta_{L_w}).
\]
Then by the fundamental sequence of class field theory (see e.g. \cite[Thm 12.1.8(iii)]{CTS-BrauerBook}, we have that, since $\beta \in \Br(L)$, then
\[
\sum_{v} \sum_{w | v} \inv_v \alpha_{k_v}\left(N_{L_w/k_v}((f^\sigma_{L_w})_\ast z_w) \right) = e
\sum_{w} \inv_w(\beta_{L_w}) = 0,\]
which proves the inclusion $\subseteq$.

For the converse inclusion, let $z = (z_v) \in Z_{0,\A}^d(X)^{\Delta_{Y/X}(\Pic(G))}$. Let $S$ be a finite set of places of $k$, containing all the infinite places, such that, over $\mathcal O_{k, S}$, we have a smooth separated model $\mathcal X$ of finite type of $X$, a smooth affine model $\mathcal G$ of $G$ with geometrically connected fibres, and  that the $G$-torsor $f: Y \to X$ extends to a $\mathcal G$-torsor $\mathcal Y \to \mathcal X$. Extend $S$ further, if necessary, so that every $x_v \in supp(z_v)$ extends to $\mathcal X(\mathcal O_{k_v(x_v)})$ for $v \notin S$, and so that the finitely many elements of $\Delta_{Y/X}(\Pic(G))$ extend to $\Br \mathcal X$.
We note that, for $v \notin S$, we then have $\alpha_{k_v}(z_v) = \sum_{x_v} n_{x_v} \cores_{k_v(x_v)/k_v}(\alpha_{k_v}(x_v)) = 0$, using the models above and since $\Br(\mathcal O_{k_v(x_v)}) = 0$. Finally, if necessary, we extend the finite set $S$ more using Lang-Weil on the special fibres with Hensel's lemma, and fix $r_v \in \mathcal X( \mathcal O_{k_v})$ for all $v \notin S$.

Now let $v \in S$ and write $z_v = z_v^+ - z_v^-$ where $z_v^+$ and $z_v^-$ are both effective, and let $\tilde{x}_v \in X_{k_v}$ be any closed point. We let
$z_v' := z_v^{-} + \tilde{x}_v$ for all $v \in S$ and let $M := (1 + |d|) lcm_{v \in S} \deg(z_v')$. Then we let
\[ {z_v^{-}}'' := \left(\frac{\# \Pic(G) \cdot M}{\deg(z_v')}\right) z_v'\]
and ${z_v^{+}}'' := z_v + {z_v^{-}}''$. (We recall that, since $G$ is connected, $\Pic(G)$ is finite.) By construction, $\left(\frac{\# \Pic(G) \cdot M}{\deg(z_v')}\right)$ is a positive multiple of $\# \Pic(G)$, and both ${z_v^{-}}''$ and ${z_v^{+}}''$ are effective, with degrees $d_- := \deg({z_v^{-}}'') = \# \Pic(G) \cdot M$ and $d_+ := \deg({z_v^{+}}'') = \# \Pic(G) \cdot M + d$ for all $v \in S$. Moreover, if $\alpha := \Delta_{Y/X}(u)$ for $u \in \Pic(G)$, then by construction $\alpha_{k_v}({z_v^{-}}'') = 0$ (since $\Delta_{Y/X}$ is additive) and $\alpha_{k_v}({z_v^{+}}'') = \alpha_{k_v}(z_v)$. Hence, by orthogonality of $z$ and using that, for $v \notin S$, the evaluation on the integral support of $z_v$ is always 0, we get
\(
\sum_{v \in S} \inv_v \alpha_{k_v}({z_v^{-}}'') = 0
\) and 
\(
\sum_{v \in S} \inv_v \alpha_{k_v}({z_v^{+}}'') = \sum_{v \in S} \inv_v \alpha_{k_v}(z_v) = 0.
\)

Now let $F/k$ be the finite Galois extension as in Lemma \ref{lem0}. We apply Lemma \ref{lem1} for $F/k$ and separately on each of the effective 0-cycles   $({z_v^{+}}'')_{v \in S}$ and  $({z_v^{-}}'')_{v \in S}$ to get finite field extensions $L_{\pm}/k$ linearly disjoint from $F$ over $k$ of degree $d_{\pm}$, and adelic points $x_\pm = (x_{\pm, w})_w \in X(\A_{L_{\pm}})$ such that, for any $v \in S$, we have
\[
\sum_{w|v} N_{L_{\pm, w}/k_v}(x_{\pm, w}) = {z^{\pm}}''_v.
\]
Hence, setting $x_{\pm, w} := (r_v)_{L_{\pm, w}}$ for $v \notin S$ and for $w|v$, we get
\[
\langle x_\pm, \alpha_L \rangle_{BM} = \sum_{v \in S} \sum_{w | v} \inv_w \alpha_{L_w}(x_{\pm, w}) = \sum_{v \in S} \inv_v \alpha_{k_v}({z^{\pm}}''_v) = 0.
\]
By the properties of $F/k$ and Lemma \ref{lem0}, as $u$ varies we have that the classes $\alpha_L$ cover the whole of $\Delta_{Y_L/X_L}(\Pic(G_L))$. Hence, by Lemma \ref{lem:lift}, we have some $\sigma_\pm \in H^1(L_\pm, G)$ and some adelic point $y_\pm \in Y^{\sigma_\pm}_{L_\pm}(\A_{L_\pm})$ such that $f^{\sigma_{\pm}}_{L_\pm}(y_\pm) = x_\pm$. 
Finally, let
\[t := z - N_{L_+/k}^\A(x_+) +  N_{L_-/k}^\A(x_-).\]
Then, for $v \in S$, we have $t_v = z_v - {z_v^{+}}'' + {z_v^{-}}'' = 0$, and for any $v \in \Omega_k$, we have $\deg(t_v) = 0$ (for $v \notin S$, we have $t_v= z_v - d r_v$). So we just need to correct the 0-cycles at $v \notin S$. For $v \notin S$, write $t_v = \sum_{s_v} n_{s_v} s_v$. The morphism $\Spec k_v(s_v) \to X$ extends to a morphism $\tilde{s}_v: \Spec \mathcal O_{k_v(s_v)} \to \mathcal X$, and by Lang's theorem \cite[Thm 2]{Lang56} and Hensel's lemma, we have $H^1(\mathcal O_{k_v(s_v)}, \mathcal G) = 0$. A section of $\tilde{s}_v^\ast \mathcal Y$ has generic fibre above $s_v$ with say image $a_v$, and $k_v(a_v) = k_v(s_v)$, so $f_\ast(a_v) = s_v$, i.e. we can lift each point $s_v$ in the support of $t_v$ to some point  $a_v \in Y_{k_v}$ with $k_v(a_v) = k_v(s_v)$. Let 
\(t'_v := \sum_{s_v} n_{s_v} a_v\)
for $v \notin S$ and $t'_v:= 0$ for $v \in S$. Then, for all $v \in \Omega_k$, we have that $\deg(t'_v) = 0$ and that $f_{k_v, \ast}(t'_v) = t_v$. Hence, the 0-cycles $y_+ \in Z_{0,\A}^{1}(Y^{\sigma_+}_{L_+})$, $-y_- \in Z_{0,\A}^{-1}(Y^{\sigma_-}_{L_-})$, and $t' \in Z_{0,\A}^{0}(Y_{k})$ recombine to $N_{L_+/k}^\Omega(x_+) -  N_{L_-/k}^\Omega(x_-) + t = z$. Hence, $z \in Z_{0,\A}^d(X)^f$, as required.
\end{proof}

\begin{corollary}\label{cor:zhang}
  Let $X$ be a smooth, quasi-projective, geometrically integral variety over a number field $k$. Let $G$ be a connected linear algebraic group over $k$ and let $f: Y \to X$ be a torsor under $G$. Then, for any $d \in \Z$, we have
  \[ Z_{0, \A}^d(X)^f =  Z_{0, \A}^d(X)^f_{ab},\]
  where $ Z_{0, \A}^d(X)^f_{ab}$ is the abelianised descent set as in \cite{Zhang-0cycles-25}.
\end{corollary}
\begin{proof}
    This follows immediately from \cite[Thm 4.3]{Zhang-0cycles-25} and from Theorem \ref{thm:main}, since Zhang proves that $ Z_{0, \A}^d(X)^f_{ab} =  Z_{0, \A}^d(X)^{\Delta_{Y/X}(\Pic(G))}$.
\end{proof}

\begin{corollary}\label{maincor}
     Let $X$ be a smooth, quasi-projective, geometrically integral variety over a number field $k$. Let $d \in \Z$. Then $Z_{0,\A}^d(X)^{\conn} = Z_{0,\A}^d(X)^{\Br}$ and $Z_{0,\A}^d(X)^{\text{\normalfont tori}} = Z_{0,\A}^d(X)^{\Br_1}$.
\end{corollary}
\begin{proof}
    By Theorem \ref{thm:main}, it suffices to show that \[\langle \Delta_{Y/X}(\Pic(G)) : [f: Y \to X] \in H^1(X, G), G \textrm{ connected linear alg group over $k$} \rangle = \Br X.\]
    In fact, using the Azumaya algebras definition of the Brauer group (via a result of Gabber, see e.g. \cite[Thms 3.3.1, 3.3.2]{CTS-BrauerBook}) and the central extensions for $n \geq 1$
    \[ 1 \to \G_{m, k} \to GL_{n, k} \to PGL_{n, k} \to 1,\]
    we have that already $\langle \Delta_{Y/X}(\Pic(PGL_{n, k})) : [f: Y \to X] \in H^1(X, PGL_{n, k}), n \geq 1 \rangle = \Br X$, showing that indeed $Z_{0,\A}^d(X)^{\conn} = Z_{0,\A}^d(X)^{\Br}$.

    To show that $Z_{0,\A}^d(X)^{\text{\normalfont tori}} = Z_{0,\A}^d(X)^{\Br_1}$, we start with the inclusion $\supseteq$.
    Any torus $G$ has $\Pic(G_{\overline{k}}) = 0$. Hence, if $u \in \Pic G$ corresponds to some central extension,  since the coboundary map $\partial_u$ is compatible with base change we have that $\partial_u(H^1(X, G)) \subset \Br_1(X):= \ker(\Br X \to \Br \overline{X})$. Hence, 
    \[\langle \Delta_{Y/X}(\Pic(G)) : [f: Y \to X] \in H^1(X, G), G \textrm{ torus over $k$} \rangle \subset  \Br_1 X\]
    and so $Z_{0,\A}^d(X)^{\text{\normalfont tori}} \supseteq Z_{0,\A}^d(X)^{\Br_1}$. To show the converse inclusion, let $\alpha \in \Br_1(X)$. Since $\alpha_{\overline{k}} = 0$, there is a finite extension $L/k$ such that $\alpha_L = 0$. Let $T := \mathfrak{R}_{L/k}(\G_{m,L})/\G_{m,k}$, where $\mathfrak{R}_{L/k}$ is the Weil restriction, so that $T$ is a torus over $k$, and consider the central extension
    \[
    1 \to \G_{m, k} \to \mathfrak{R}_{L/k}(\G_{m,L}) \to T \to 1,
    \]
    which can be viewed as an element $u \in \Pic(T)$.
    By Shapiro's lemma, we have that $H^i(X, \mathfrak{R}_{L/k}(\G_{m,L})) = H^i(X_L, \G_{m,L})$, and the map $H^2(X, \G_m) \to H^2(X, \mathfrak{R}_{L/k}(\G_{m,L}))$ coming from the extension above is the restriction map $\res_{L/k}: \Br X \to \Br(X_L)$. In fact, we have the exact sequence 
    \[ \Pic(X_L) \to H^1(X, T) \xrightarrow{\partial_u} \Br(X) \xrightarrow{\res_{L/k}} \Br(X_L).\]
    Since $\alpha_L = \res_{L/k}(\alpha) = 0$, we have by exactness that $\alpha = \partial_u(\tau)$ for some $\tau \in H^1(X, T)$. Hence, $\alpha \in \langle \Delta_{Y/X}(\Pic(G)) : [f: Y \to X] \in H^1(X, G), G \textrm{ torus $k$} \rangle$. Hence, $Z_{0,\A}^d(X)^{\text{\normalfont tori}} \subseteq Z_{0,\A}^d(X)^{\Br_1}$, which proves the other inclusion and thus the result, as required.
\end{proof}

\begin{remark}
    In fact, a similar proof also shows that, for rational points, $X(\A_k)^{\text{tori}} = X(\A_k)^{\Br_1}$. 
\end{remark}

\section{A reformulation of Conjecture (E) using connected descent}
In this section, we take $X$ to be projective. In particular, $Z^d_{0, \Omega}(X) = Z^d_{0, \A}(X)$ for all $d \in \Z$. We also note that if $f: Y \to X$ is a torsor under a connected linear algebraic group $G$ over $k$, then 
\[ Z^d_{0, \A}(X)^f \subset Z^d_{0, \Omega}(X)^f \subset Z^d_{0, \Omega}(X)^{\Delta_{Y/X}(\Pic(G))} = Z^d_{0, \A}(X)^{\Delta_{Y/X}(\Pic(G))}.\]
In the light of Theorem \ref{thm:main}, it is natural to ask if there is a formulation of Conjecture (E) in terms of the connected descent for 0-cycles. We prove the following, answering in the positive.

\begin{theorem}\label{conjE} Let $X$ be a smooth, projective, geometrically integral variety over a number field $k$. 
Let $\phi: Z_{0, \Omega}(X) \to \widehat{\prod_v \CH'_0(X_v)}$ be given by $(z_v)_v \mapsto (([z_{v}]_n)_{n \geq 1})_v$, where $[- ]_n$ denotes the class in $\CH'_0(X_{k_v})/n$.
    Then Conjecture (E) is equivalent to 
    \[
    \im\left(\widehat{\CH_0(X)} \to  \widehat{\prod_v \CH'_0(X_v)} \right) = \bigcap_{G \textrm{ connected lin alg $k$-groups}} \bigcap_{[f: Y \to X] \in H^1(X, G)} \overline{\phi(Z_{0, \Omega}(X)^f)}.
    \]
    In particular, it suffices to consider $G = PGL_{n, k}$ for $n \geq 1$. 
\end{theorem}

Before we can prove Theorem \ref{conjE}, we need some preliminary results.

\begin{proposition} \label{propconje} Let $X$ be a smooth, projective, geometrically integral variety over a number field $k$. Let $G$ be a connected linear algebraic group over $k$ and let $f: Y \to X$ be a torsor under $G$. Let $u \in \widehat{\prod_v \CH'_0(X_v)}$. Then the following are equivalent:
\begin{enumerate}
    \item $u \in (\Br_0(X) + \Delta_{Y/X}(\Pic(G)))^{\perp}$, where $(-)^\perp$ denotes the annihilator in $\widehat{\prod_v \CH'_0(X_v)}$.
    \item For any $n \geq 1$, there is a $d \in \Z$ and $z \in Z_{0, \Omega}^d(X)^f$ such that $[z_{v}]_n= u_{v, n}$ for all $v \in \Omega_k$, where $[-]$ denotes the class in $\CH'_0(X_{k_v})$ and the subscript $n$ denotes the image in $\CH'_0(X_{k_v})/n$.
    \item $u \in \overline{\phi(Z_{0, \Omega}(X)^f)}$, where $\phi: Z_{0, \Omega}(X) \to \widehat{\prod_v \CH'_0(X_v)}$.
\end{enumerate}
    
\end{proposition}

\begin{proof} $(1) \Rightarrow (2)$: We fix a finite place $v_0$. At any finite place $v$, the degree map induces $\deg(u_v) \in \widehat{\Z}$. Let $v \neq v_0$ be finite and let $m \geq 2$. By the fundamental sequence of class field theory, there is some $a \in \Br k$ such that $\inv_{v} a = 1/m$, $\inv_{v_0} a = -1/m$, and $\inv_{v'} a = 0$ for all $v' \neq v, v_0$. Since by assumption $u \in (\Br_0(X) + \Delta_{Y/X}(\Pic(G)))^{\perp}$, by orthogonality we have that
\(0 = \langle u, a_X \rangle_{BM} = \frac{1}{m}(\deg(u_v) - \deg(u_{v_0})) \in \Q/\Z.\)
Since this holds for any $m \geq 2$, it follows that $\deg(u_v) = \deg(u_{v_0}) =: \delta$ for all finite places $v$. For a real place $v$, a similar argument with $\inv_v(a) = 1/2$, $\inv_{v_0}(a) = -1/2$ and $\inv_{v'} (a) = 0$ for all $v' \neq v_0, v$ shows that $\deg(u_v) = \delta \bmod 2$. 

Now let $n \geq 1$. Since $\# \Delta_{Y/X}(\Pic(G))$ is finite, we can take $N:= lcm(n, \# \Delta_{Y/X}(\Pic(G)))$. 
Let $\tilde{x} \in X$ be a closed point with $\tilde{d} := [k(\tilde{x}): k]$, and let $d \in \Z$ be such that $d \equiv \delta \bmod (2 N \tilde{d})$. Let $v$ be a finite place and let $x_v \in Z_0(X_{k_v})$ be a 0-cycle representing $u_{v, N \tilde{d}} \in \CH_0'(X_{k_v})/(N \tilde{d})$. Then $N \tilde{d} | d - \deg(x_v)$. Let
\[
z_v := x_v + N \left(\frac{d - \deg(x_v)}{N \tilde{d}} \right) \tilde{x}_{k_v}.
\]
Then $\deg z_v = d$ and $[z_v]_N = u_{v, N}$. 

Now let $v$ be a real place. We note that, in this case,  $\CH'_0(X_v)$ is finite and 2-torsion (see e.g. \cite[\S2.1]{Liang-Arithmetic-13}), and hence $\CH'_0(X_v) = \widehat{\CH'_0(X_v)}$. Let  $x_v \in Z_0(X_{k_v})$ be a 0-cycle representing $u_{v}$, and let $x'$ be a closed point in $ X_{\overline{k}_v}$ where $\overline{k}_v \cong \C$. Then $2 | d - \deg(x_v)$ and we let
\[
z_v := x_v + \left(\frac{d - \deg(x_v)}{2} \right) N_{\overline{k}_v/k_v}(x').
\]
Then $\deg z_v = d$ and $[z_v] = u_v$.  At a complex place $w$, we choose $x'_w \in X(k_w)$ and take $z_w = d x'_w$, since $\CH'_0(X_{k_w}) = 0$.
Hence, we have constructed $z := (z_v)_v \in Z_{0, \Omega}^d(X)$ such that $[z_v]_N = u_{v, N}$ for all places $v \in \Omega_k$. Moreover, for any $\alpha \in \Delta_{Y/X}(\Pic(G))$, the evaluation at $\alpha$ factors through the local Chow groups modulo the order of $\alpha$, and $ord(\alpha) | N$. Hence, 
\(\langle z, \alpha \rangle_{BM} = \langle u, \alpha \rangle_{BM}
\)
for all $\alpha \in \Delta_{Y/X}(\Pic(G))$ and, by Theorem \ref{thm:main}, we have that $z \in Z^d_{0, \Omega}(X)^f$. Moreover, since $n | N$ by definition, we can reduce from modulo $N$ to modulo $n$ to get $[z_v]_n = u_{v, n}$ for all $v \in \Omega_k$, thus proving (2). \\

\noindent $(2) \Rightarrow (3)$: Since $u \in \widehat{\prod_v \CH'_0(X_v)}$, a basis of neighbourhoods for $u$ in $\widehat{\prod_v \CH'_0(X_v)}$ is given by
\[
U(u; S, (n_v)_{v\in S}):= u + \bigcap_{v \in S} \ker\left(q_{v, n_v}: \prod_{v \in \Omega_k}\widehat{ \CH'_0(X_v)} \to \CH'_0(X_v)/n_v\right)
\]
ranging over all finite subsets of places $S \subset \Omega_k$, and all tuples of positive integers $(n_v)_{v \in S}$. Fix such a basic neighbourhood for some $S$ and some $(n_v)_{v \in S}$, and let $n := lcm_{v \in S} n_v$. Since we are assuming (2), we have some $d \in \Z$ and some $z \in Z_{0, \Omega}^d(X)^f$ such that $[z_v]_n = u_{v, n}$ for all $v \in \Omega_k$. For each $v \in S$, we let $\pi_{v; n, n_v}:\CH'_0(X_v)/n \to \CH'_0(X_v)/n_v$ be the natural reduction map. Then
\[
q_{v, n_v}(\phi(z) - u) = \pi_{v; n, n_v}([z_v]_n - u_{v, n}) = 0
\]
for each $v \in S$. Hence, $\phi(z) \in U(u;S,(n_v)_{v \in S}) \cap \phi(Z_{0, \Omega}(X)^f)$. Since this was true for an arbitrary basic neighbourhood of $u$, we have that $u \in \overline{\phi(Z_{0, \Omega}(X)^f)}$, which proves (3).\\

\noindent $(3) \Rightarrow (1)$: Let $\alpha \in \Br_0(X) + \Delta_{Y/X}(\Pic(G))$ and let $n:= ord(\alpha)$ be its order. Let $S$ be a finite set of places of $k$, including the infinite places, such that $\alpha$ extends to an element in $\Br(\mathcal X)$ for some proper model $\mathcal X$ of $X$ over $\mathcal O_{k, S}$. Then, for any $v \notin S$ and any closed point $x_v \in X_{k_v}$, by properness and using the fact that $\Br(\mathcal O_{k_v(x_v)}) = 0$, we have that the evaluation $\alpha(x_v) = 0$.
Now let $v \in S$. We note that the local pairing on $Z_{0}(X_{k_v})$ is compatible with descending to $\CH'_0(X_{k_v})$ (see e.g. \cite[\S14.1]{CTS-BrauerBook}). Moreover, for any $t \in \CH'_0(X_{k_v})$, we have that $\langle n t, \alpha\rangle_{v} = \langle  t, n \alpha\rangle_{v} = 0$, so the local pairing only depends on the image of $t$ in $\CH'_0(X_{k_v})/n$.

Assuming (3), let $u \in \overline{\phi(Z_{0, \Omega}(X)^f)}$. Consider the open basic neighbourhood of $u$ given by $U(u; S, (n)_{v \in S})$. By definition of the closure, every open neighbourhood of $u$ meets $\phi(Z_{0, \Omega}(X)^f)$, and hence in particular
\[
U(u; S, (n)_{v \in S}) \cap \phi(Z_{0, \Omega}(X)^f) \neq \emptyset.
\]
This implies that there is some $z \in Z_{0, \Omega}(X)^f$ such that $\phi(z) \in U(u; S, (n)_{v \in S})$. By definition of  $U(u; S, (n)_{v \in S})$, we have that 
\[
0 = q_{v, n}(\phi(z) - u) = q_{v, n}(\phi(z)) - q_{v, n}(u) = [z_v]_n - u_{v, n}
\]
for any $v \in S$. Hence, $z$ and $u$ have the same local pairing with $\alpha$ for all $v \in S$, since $[z_v]_n = u_{v, n}$ and the local pairing only depends on the image in $\CH'_0(X_{k_v})/n$. This implies that
\[
\langle u, \alpha \rangle_{BM} = \sum_{v \in S} \langle u_{v, n}, \alpha \rangle_v = \sum_{v \in S} \langle [z_v]_n, \alpha \rangle_v = \langle z, \alpha \rangle_{BM} = 0,
\]
where the right-most equality follows since, by  Theorem \ref{thm:main},  any $z \in Z^d_{0, \Omega}(X)^f$ annihilates $\Delta_{Y/X}(\Pic(G))$ and, moreover, it also annihilates $\Br_0(X)$ since $\deg(z_v) = d$ for all $v \in \Omega_k$ and by the fundamental exact sequence of class field theory. Since $\alpha$ was arbitrary, we have that $u \in (\Br_0(X) + \Delta_{Y/X}(\Pic(G)))^\perp$, which proves (1), as required.
\end{proof}

\begin{proof}[Proof of Theorem \ref{conjE}]
This follows immediately from Proposition \ref{propconje} by taking the intersections over all connected linear groups $G$ and all torsors $f: Y \to X$ under $G$
\[
\begin{array}{ll}
\bigcap_{G \text{ conn}} \bigcap_{[f] \in H^1(X, G)} \overline{\phi(Z_{0, \Omega}(X)^f)} & = \bigcap_{G \text{ conn}} \bigcap_{[f] \in H^1(X, G)} (\Br_0(X) + \Delta_{Y/X}(\Pic(G)))^\perp\\
& = \left(\Br_0(X) + \sum_{G, [f]} \Delta_{Y/X}(\Pic(G))\right)^\perp\\ 
&= (\Br X)^\perp, 
\end{array}
\]
where we have used the fact that the $\Delta_{Y/X}(\Pic(G))$ generate $\Br X$ (see the proof of Corollary \ref{maincor}). In fact, taking $\Delta_{Y/X}(\Pic(G))$ ranging over $G = PGL_{n, k}$ for $n \geq 1$ suffices. The equivalence with Conjecture (E) follows immediately since Conjecture (E) is the statement that
\[
    \im\left(\widehat{\CH_0(X)} \to  \widehat{\prod_v \CH'_0(X_v)} \right) = (\Br X)^\perp,
    \]
as required.
\end{proof}

\printbibliography

\end{document}